\documentclass[11pt,dvips]{article}
\usepackage{amsmath,amsfonts,amssymb,amsthm,graphicx,verbatim,subfig}
\usepackage[all]{xy}
\usepackage[a4paper,left=25mm]{geometry}
\ifx\pdftexversion\undefined

\usepackage[a4paper,colorlinks,link
=black,filecolor=black,citecolor=black,urlcolor=black,pdfstartview=FitH]{hyperref}
\else

\usepackage[a4paper,colorlinks,linkcolor=black,filecolor=black,citecolor=black,urlcolor=black,pdfstartview=FitH]{hyperref}
\fi

\font\sixbb=msbm6
\font\eightbb=msbm8
\font\twelvebb=msbm10 scaled 1095
\newfam\bbfam
\textfont\bbfam=\twelvebb \scriptfont\bbfam=\eightbb
                           \scriptscriptfont\bbfam=\sixbb
\def\bb{\fam\bbfam\twelvebb}

\newcommand{\FF}{{\bb F}}

\newtheorem{theorem}{\bf Theorem}[section]
\newtheorem{claim}[theorem]{\bf Claim}

\newtheorem{proposition}[theorem]{\bf Proposition}
\newtheorem{corollary}[theorem]{\bf Corollary}
\newcommand{\enp}{\begin{flushright} $\Box$ \end{flushright}}
\newcommand{\beq}[0]{\begin{equation}}
\newcommand{\enq}[0]{\end{equation}}

\newcommand{\cb}{{\cal B}}
\newcommand{\cm}{{\cal M}}
\newcommand{\rk}{{\rm rank}}

\newcommand{\ccs}{\mathcal{S}}
\newcommand{\ccb}{\mathcal{B}}
\newcommand{\ccm}{\mathcal{M}}
\newcommand{\ccp}{\mathcal{P}}
\newcommand{\ccq}{\mathcal{Q}}
\newcommand{\pf}{\text{Pf}}

\newcommand{\csy}{{\mathbb S}}
\newcommand{\nle}{[n]_{\leq}^2}
\newcommand{\qtil}{\tilde{q}}
\newcommand{\ktil}{\tilde{K}_n}

\newcommand{\sgn}{\text{sgn}}

\title{Maximal Rank in Matrix Spaces \\ via Graph Matchings}
\begin{document}
\author{Roy Meshulam\thanks{Department of Mathematics,
Technion, Haifa 32000, Israel. e-mail:
meshulam@math.technion.ac.il~. Supported by ISF grant no. 326/16 and GIF grant no. 1261/14.}}
\maketitle
\pagestyle{plain}

\begin{abstract}
Let $M_{n}(\FF)$ be the space of $n \times n$ matrices over a field $\FF$.
A matrix $A=\big(A(i,j)\big)_{i,j=1}^n \in M_n(\FF)$ is weakly symmetric if $A(i,j) \neq 0$ iff $A(j,i) \neq 0 $ holds for all $i,j$. A matrix is alternating if it is skew-symmetric with zero diagonal. Let $W_n(\FF)$ and $A_n(\FF)$ denote respectively the set of  weakly symmetric matrices and the space of alternating matrices in $M_n(\FF)$.
Let $[n]=\{1,\ldots,n\}$. For $0 \neq A \in W_n(\FF)$ let $\qtil(A)=\{i,j\}$, where $(i,j)$ is the unique pair in $[n]^2$
such that $A(i,j) \neq 0$ and $A(i',j')=0$ whenever $j<j'$ or $j=j'$ and $i<i'$.
For a translate $\ccs$ of a linear space $\cb \subset W_n(\FF)$ let
$G_{\ccs}$ be the graph with loops on the vertex set $[n]$  with edge set $E_{\ccs}=\{\qtil(B):0 \neq B \in \cb\}$.
A subset $M \subset E_{\ccs}$ is a matching if $e \cap e'=\emptyset$ for all $e \neq e' \in M$.
Let $\mu(G_{\ccs})=\max\sum_{e \in M} |e|$ where $M$ ranges over all matchings $M \subset E_{\ccs}$. Let $\rho(\ccs)$ denote the maximal rank of a matrix in $\ccs$.

It is shown that if $\ccs$ is a translate of a linear space contained in $W_n(\FF)$ and $|\FF|\geq 3$ then
$\rho(\ccs)\geq \mu(G_{\ccs})$. The restriction on $\FF$ can be removed if $\ccs$ is an affine subspace of $A_n(\FF)$.
Applications include simple proofs of upper bounds on the dimension of affine subspaces of symmetric and alternating matrices of bounded rank.

\end{abstract}

\section{Introduction}
\label{s:intro}

Let $M_{n}(\FF)$ be the space of $n \times n$ matrices over a field $\FF$.
For a subset $\ccs \subset M_n(\FF)$ let $\rho(\ccs)=\max\{\rk(A):A \in \ccs\}$ denote the maximum rank of a matrix in $\ccs$.
Let $H_n(\FF)$ denote the space of symmetric matrices in $M_n(\FF)$.  A matrix $A=\big(A(i,j)\big)_{i,j=1}^n \in M_n(\FF)$
is {\it alternating} if $A=-A^T$ and $A(i,i)=0$ for $1 \leq i \leq n$. Let $A_n(\FF)$ denote the space of alternating matrices in $M_n(\FF)$.
A matrix $A \in M_n(\FF)$ is
{\it weakly symmetric} if $A(i,j) \neq 0$ iff $A(j,i) \neq 0 $ holds for all $i,j$. Let $W_n(\FF)$ denote the set of all weakly symmetric matrices in $M_n(\FF)$. Note that $A_n(\FF), H_n(\FF) \subset W_n(\FF)$ and
$W_n(\FF_2)=H_n(\FF_2)$. In this note we study lower bounds on $\rho(\ccs)$ for affine translates $\ccs$ of linear spaces of weakly symmetric matrices, in terms of matching numbers of a certain graph associated with $\ccs$.

Let $[n]=\{1,\ldots,n\}$ and let $\nle=\{(i,j) \in [n]^2: i \leq j\}$.  The {\it colexicographic order} on
$\nle$ is given by $(i,j) \prec (i',j')$ iff $j <j'$ or $j=j'$ and $i<i'$. Equivalently,
$(i,j) \prec (i',j')$ iff $2^i+2^j <2^{i'}+2^{j'}$.
Let $K_n=\{e \subset [n]:|e|=2\}$ denote the edge set of the complete graph on $[n]$ and let
$\ktil=\{e \subset [n]: 0<|e| \leq 2\}$ denote the edge set of the complete graph with loops on $[n]$.
 A subset $M  \subset \ktil$ is a {\it matching} if $e \cap e' =\emptyset$ for all $e \neq e' \in M$.
For a graph with loops $G \subset \ktil$ let $\cm(G)$ denote the set of all matchings $M \subset G$.
Let $\nu(G)=\max \{|M|: M \in \cm(G)\}$ and
$\mu(G)=\max \{\sum_{e \in M} |e|: M \in \cm(G)\}$. A matching $M \subset \ktil$ is {\it perfect} if $\mu(M)=n$. Note that if $G$ is loopless, i.e. $G \subset K_n$, then $\nu(G)$ is the usual matching number of $G$ and $\mu(G)=2\nu(G)$.

For $0 \neq A=\big(A(i,j)\big)_{i,j=1}^n \in W_n(\FF)$ let
$q(A)= \max\{(i,j) \in \nle:A(i,j)\neq 0\}$ where the maximum is taken with respect to the colexicographic order.
For $A$ such that $q(A)=(i,j)$ let $\qtil(A)=\{i,j\} \in \ktil$.
Let $\ccs$ be a translate of a linear space $\ccb \subset W_n(\FF)$, i.e. $\ccs=A+\ccb$ for some $A \in M_n(\FF)$.
Associate with $\ccs$ a graph with loops
$$G_{\ccs}=\{\qtil(B): 0 \neq B \in \ccb\}=\{\qtil (S_1-S_2):S_1\neq S_2 \in \ccs\} \subset \ktil.$$
Our main results provide a link between the maximum rank in $\ccs$ and matchings in $G_{\ccs}$.
\begin{theorem}
\label{t:rankmatch}
Suppose $|\FF| \geq 3$ and let $\ccs=A+\ccb$ where $A \in M_n(\FF)$ and $\ccb$ is a linear space contained in $W_n(\FF)$.
Then $\rho(\ccs) \geq \mu(G_{\ccs})$.
\end{theorem}
\noindent
The restriction on $\FF$ is superfluous if $\ccs$ is an affine space of alternating matrices:
\begin{theorem}
\label{t:alter}
Let $\FF$ be an arbitrary field and let $\ccs$ be an affine subspace of $A_n(\FF)$.
Then $\rho(\ccs) \geq \mu(G_{\ccs})$.
\end{theorem}
\noindent
{\bf Remarks:}
\\
1. The case $A=0$ and $|\FF| \geq \mu(G_{\ccb})+1$ of Theorem \ref{t:rankmatch} had (essentially) been proved in \cite{M89}.
The approach to the present improved result is somewhat different and uses additional ideas.
The proof of Theorem \ref{t:alter} utilizes Pfaffians of alternating matrices.
\\
2. Theorem \ref{t:rankmatch} does not hold for $\FF=\FF_2$ as the following examples show.
\\
(i) The affine subspace
$$
\ccs=\Big\{ \left[
\begin{array}{ll}
x & x \\
x & 1
\end{array}
\right] : x \in \FF_2\Big\} \subset H_2(\FF_2)
$$
satisfies $\mu(G_{\ccs})=2>1=\rho(\ccs)$.
\\
(ii) The linear subspace
$$
\ccs=\Big\{ \left[
\begin{array}{lll}
x & x & 0 \\
x & 0 & y \\
0 & y & y
\end{array}
\right] : x,y \in \FF_2\Big\} \subset H_3(\FF_2)
$$
satisfies $\mu(G_{\ccs})=3>2=\rho(\ccs)$.
\\
Note however that $\rho(\ccs)=\nu(G_{\ccs})$ holds in both examples. In fact, the following consequence of Theorem \ref{t:alter} holds.
\begin{corollary}
\label{c:symm}
Let $\ccs$ be an affine subspace of $H_n(\FF_2)$.
Then $\rho(\ccs) \geq \nu(G_{\ccs})$.
\end{corollary}
\noindent
Theorems \ref{t:rankmatch} and \ref{t:alter} can be used to provide short combinatorial proofs of upper bounds on the maximal dimensions of affine subspaces of $A_n(\FF)$ and of translates of linear spaces contained in $W_n(\FF)$, that consist of matrices of bounded rank. We first consider the alternating case. For $k=2t$ even and $n \geq k$ let
$$U_1^a(n,k)= \big\{\big(A(i,j)\big)_{i,j=1}^n \in A_n(\FF): A(i,j)=0 \text{~if~} \max\{i,j\}>k+1 \big\}$$
and
$$U_2^a(n,k)=\big\{\big(A(i,j)\big)_{i,j=1}^n \in A_n(\FF): A(i,j)=0 \text{~if~} \min\{i,j\}>\frac{k}{2} \big\}.$$
Then $\rho(U_1^a(n,k))=\rho(U_2^a(n,k))=k$. Define
\begin{equation*}
\label{e:dua}
\begin{split}
u_a(n,k)&= \max \big\{\dim U_1^a(n,k),\dim U_2^a(n,k)\big\} \\
&= \max \Big\{\binom{2t+1}{2},tn-\binom{t+1}{2}\Big\}.
\end{split}
\end{equation*}
\noindent
The following result was first obtained in \cite{M89} under the assumption that $\ccs$ is a linear subspace of $A_n(\FF)$ and $|\FF|\geq n+1$. The general case has recently been proved by de Seguins Pazzis \cite{DSP16}. Here it is obtained as a direct consequence of Theorem \ref{t:alter}.
\begin{theorem}[\cite{M89,DSP16}]
\label{t:alt1}
Let $\FF$ be any field and let $\ccs$ be an affine subspace of $A_n(\FF)$ such that
$\rho(\ccs)=k<n$. Then $\dim \ccs \leq u_a(n,k)$.
\end{theorem}
\noindent
We next consider the weakly symmetric case.
Let $k \leq n$. Define
$$U_1^s(n,k)= \big\{\big(A(i,j)\big)_{i,j=1}^n \in H_n(\FF): A(i,j)=0 \text{~if~} \max\{i,j\}>k \big\}.$$
For $k=2t$ even let
$$
U_2^s(n,k)=
\big\{\big(A(i,j)\big)_{i,j=1}^n \in H_n(\FF): A(i,j)=0 \text{~if~} \min\{i,j\}>t \big\},
$$
and for $k=2t+1$ odd let
$$
U_2^s(n,k)=
\big\{\big(A(i,j)\big)_{i,j=1}^n \in H_n(\FF): A(i,j)=0 \text{~if~} \min\{i,j\}>t ~\&~ (i,j)\neq (t+1,t+1) \big\}.
$$
Then $\rho(U_1^s(n,k))=\rho(U_2^s(n,k))=k$. Define
\begin{equation*}
\label{e:dus}
\begin{split}
u_s(n,k)&=\max \big\{\dim U_1^s(n,k),\dim U_2^s(n,k)\big\} \\
&= \left\{
\begin{array}{ll}
\max\big\{\binom{2t+1}{2},tn-\binom{t}{2}\big\} & k=2t, \\
\max\big\{\binom{2t+2}{2},tn-\binom{t}{2}+1\big\} & k=2t+1.
\end{array}
\right.
\end{split}
\end{equation*}
\noindent
The following result was first obtained in \cite{M89} under the assumption that $\ccs$ is a linear subspace of $H_n(\FF)$ and
$|\FF|\geq n+1$. It has recently been proved by de Seguins Pazzis \cite{DSP16} for affine subspaces $\ccs$ of $H_n(\FF)$ with no restrictions on $\FF$. Here it is established for slightly more general affine spaces provided that $\FF \neq \FF_2$.
\begin{theorem}
\label{t:ws1}
Suppose $|\FF|\geq 3$ and
let $\ccs=A+\cb$ where $A \in M_n(\FF)$ and $\cb$ is a linear space contained in $W_n(\FF)$.
If $\rho(\ccs)=k$ then $\dim \ccs \leq u_s(n,k)$.
\end{theorem}
\ \\ \\
The paper is organized as follows. Theorem \ref{t:rankmatch} is proved in Section \ref{s:wsym}.
In Section \ref{s:alter} we use Pfaffians to establish Theorem \ref{t:alter} and then deduce Corollary \ref{c:symm}.
In Section \ref{s:eg} we recall the Erd\H{o}s-Gallai theorem \cite{EG61} on the maximal number of edges in simple graphs with bounded matching number, as well as the analogous  result for graphs with loops \cite{M89}. Combining these with Theorems
\ref{t:alter} and \ref{t:rankmatch}, directly implies Theorems \ref{t:alt1} and \ref{t:ws1}.
We conclude in Section \ref{s:conc} with some remarks and open problems.

\section{Translates of Weakly Symmetric Subspaces}
\label{s:wsym}

In this section we prove Theorem \ref{t:rankmatch}. Let $\ccs=A+\ccb$ be an affine subspace of $M_n(\FF)$ where $A \in M_n(\FF)$ and $\ccb \subset W_n(\FF)$ is a linear space.
We first note that by restricting to the principal submatrix determined by a matching $M \subset G_{\ccs}$ such that
$\sum_{e \in M}|e|=\mu(G_{\ccs})$, it suffices to show that if $\mu(G_{\ccs})=n$, then $\ccs$ contains a nonsingular matrix.
Choose $B_1,\ldots,B_t \in \ccb$ such that
$M=\{\qtil(B_1),\ldots,\qtil(B_t)\}$ is a prefect matching of $\ktil$, and let $|\qtil(B_i)|=\delta_i$.
Then
$\sum_{i=1}^t \delta_i=n$. For $j=1,2$ let $I_j=\{1 \leq i \leq t: \delta_i=j\}$.
For $1 \leq i \leq t$ let
$$
q(B_i)=\left\{
\begin{array}{cc}
(\alpha_i,\alpha_i)  & i \in I_1, \\
(\beta_i,\gamma_i) & i \in I_2.
\end{array}
\right.
$$
Let $x=(x_1,\ldots,x_t)$ be a vector of variables and let
$$f(x)=\det\Big(A+\sum_{i=1}^t x_i B_i\Big) \in \FF[x_1,\ldots,x_t].$$
We have to show that there exists an $\lambda=(\lambda_1,\ldots,\lambda_t) \in \FF^t$ such that $f(\lambda) \neq 0$.
The main ingredient in the proof is the following
\begin{proposition}
\label{p:coeff}
The monomial $\prod_{i=1}^t x_i^{\delta_i}$ appears in $f(x)$ with a nonzero coefficient.
\end{proposition}
\noindent
{\bf Proof:} Let $\csy_n$ be the symmetric group on $[n]$ and let $\ccp$ denote the set of all ordered partitions $P=(C_0,C_1,\ldots,C_t)$ of $[n]$ into $t+1$ disjoint sets. Then:
\begin{equation}
\label{e:exp}
\begin{split}
f(x)&=\det(A+\sum_{i=1}^t x_i B_i) \\
&=\sum_{\sigma \in \csy_n} \sgn(\sigma) \prod_{j=1}^n \big(A(j,\sigma(j))+\sum_{i=1}^t x_i B_i(j,\sigma(j))\big) \\
&=\sum_{\sigma \in \csy_n} \sgn(\sigma) \sum_{P=(C_0,\ldots,C_t) \in \ccp} \big( \prod_{j \in C_0} A(j,\sigma(j))
\prod_{i=1}^t  \prod_{j \in C_i} B_i(j,\sigma(j))\big) \prod_{i=1}^t x_i^{|C_i|}.
\end{split}
\end{equation}
Let $\ccq$ denote the set of all ordered partitions $Q=(C_1,\ldots,C_t)$ of $[n]$ into $t$ disjoint sets
such that $|C_i|=\delta_i$. For $Q=(C_1,\ldots,C_t) \in \ccq$ and $\sigma \in \csy_n$ let
$$
g(Q,\sigma)=\prod_{i=1}^t  \prod_{j \in C_i} B_i(j,\sigma(j)).
$$
As $\sum_{i=1}^t \delta_i=n$,
it follows from (\ref{e:exp}) that the coefficient of $\prod_{i=1}^t x_i^{\delta_i}$ in $f(x)$ is
\begin{equation}
\label{e:exp1}
\sum_{\sigma \in \csy_n} \sgn(\sigma) \sum_{Q=(C_1,\ldots,C_t) \in \ccq}
g(Q,\sigma).
\end{equation}
Consider the partition $Q_0=(\qtil(B_1),\ldots,\qtil(B_t)) \in \ccq$ and the permutation
$\sigma_0 \in \csy_n$ given by
$$\sigma_0(j)=\left\{
\begin{array}{cc}
\alpha_i & j=\alpha_i,  \\
\gamma_i & j=\beta_i, \\
\beta_i & j=\gamma_i.
\end{array}
\right.
$$
Note that $g(Q_0,\sigma_0)\neq 0$. Hence,
Proposition \ref{p:coeff} will follow from (\ref{e:exp1}) and the next observation.
\begin{claim}
\label{c:nzero}
Let $(Q,\sigma) \in \ccq\times \csy_n$ such that $g(Q,\sigma)\neq 0$. Then
$(Q,\sigma)=(Q_0,\sigma_0)$.
\end{claim}
\noindent
{\bf Proof:}
Write $Q=(C_1,\ldots,C_t)$ where
$$
C_i=\left\{
\begin{array}{cc}
\{k_i\} & i \in I_1, \\
\{\ell_i,m_i\} & i \in I_2.
\end{array}
\right.
$$
Then
\begin{equation*}
\label{e:gqs}
\begin{split}
g(Q,\sigma)&=\prod_{i=1}^t  \prod_{j \in C_i} B_i(j,\sigma(j)) \\
&= \prod_{i \in I_1} B_i(k_i,\sigma(k_i)) \prod_{i \in I_2} \Big( B_i(\ell_i,\sigma(\ell_i))\cdot B_i(m_i,\sigma(m_i)) \Big).
\end{split}
\end{equation*}
As $g(Q,\sigma) \neq 0$ it follows that for $i \in I_1$:
\begin{equation*}
\big(\min(k_i,\sigma(k_i)),\max(k_i,\sigma(k_i)\big) \preceq q(B_i)=(\alpha_i,\alpha_i),
\end{equation*}
hence
\begin{equation}
\label{e:ki1}
2^{k_i}+2^{\sigma(k_i)} \leq 2^{\alpha_i}+2^{\alpha_i}.
\end{equation}
Similarly, for $i \in I_2$:
\begin{equation*}
\begin{split}
&\Big(\min\big(\ell_i,\sigma(\ell_i)\big),\max\big(\ell_i,\sigma(\ell_i)\big)\Big) \preceq q(B_i)=(\beta_i,\gamma_i), \\
&\Big(\min\big(m_i,\sigma(m_i)\big),\max\big(m_i,\sigma(m_i)\big)\Big) \preceq q(B_i)=(\beta_i,\gamma_i),
\end{split}
\end{equation*}
hence
\begin{equation}
\label{e:klmi1}
\begin{split}
&2^{\ell_i}+2^{\sigma(\ell_i)} \leq 2^{\beta_i}+2^{\gamma_i}, \\
&2^{m_i}+2^{\sigma(m_i)} \leq 2^{\beta_i}+2^{\gamma_i}.
\end{split}
\end{equation}
Next note that $Q_0, Q$ and $\sigma(Q)=(\sigma(C_1),\ldots,\sigma(C_t))$ are all partitions of $[n]$, therefore
\begin{equation}
\label{e:esums}
\begin{split}
&\sum_{i \in I_1} 2^{\sigma(k_i)} + \sum_{i \in I_2}  \big(2^{\sigma(\ell_i)}+2^{\sigma(m_i)}\big)\\
&=\sum_{i \in I_1} 2^{k_i} + \sum_{i \in I_2}  \big(2^{\ell_i}+2^{m_i}\big) \\
&= \sum_{i \in I_1} 2^{\alpha_i} + \sum_{i \in I_2}  \big(2^{\beta_i}+2^{\gamma_i}\big)\\
&=\sum_{j=1}^n 2^j=2^{n+1}-2
\end{split}
\end{equation}
Summing (\ref{e:ki1}) over $i \in I_1$ and (\ref{e:klmi1}) over $i \in I_2$ it follows by (\ref{e:esums}) that
\begin{equation}
\label{e:sumklm}
\begin{split}
2^{n+2}-4&=\sum_{i \in I_1} \big(2^{k_i}+2^{\sigma(k_i)} \big)+\sum_{i \in I_2}
\Big( \big(2^{\ell_i}+2^{\sigma(\ell_i)}\big)+\big(2^{m_i}+2^{\sigma(m_i)} \big) \Big) \\
&\leq \sum_{i \in I_1} \big(2^{\alpha_i}+2^{\alpha_i}\big)+
\sum_{i \in I_2} \Big(\big(2^{\beta_i}+2^{\gamma_i}\big)+\big(2^{\beta_i}+2^{\gamma_i}\big)\Big) \\
&=2^{n+2}-4.
\end{split}
\end{equation}
By (\ref{e:sumklm}), all inequalities in (\ref{e:ki1}) and (\ref{e:klmi1}) are in fact equalities. Thus, for all $i \in I_1$
\begin{equation}
\label{e:ki1f}
2^{k_i}+2^{\sigma(k_i)} = 2^{\alpha_i}+2^{\alpha_i},
\end{equation}
and for all $i \in I_2$
\begin{equation}
\label{e:klmi1f}
\begin{split}
2^{\ell_i}+2^{\sigma(\ell_i)}&=2^{\beta_i}+2^{\gamma_i}, \\
2^{m_i}+2^{\sigma(m_i)}&= 2^{\beta_i}+2^{\gamma_i}.
\end{split}
\end{equation}
It follows by (\ref{e:ki1f}) that $k_i=\sigma(k_i)=\alpha_i$ for all $i \in I_1$. Similarly,
(\ref{e:klmi1f}) implies that for all $i \in I_2$
$$\{\ell_i,\sigma(\ell_i)\}=\{\beta_i,\gamma_i\}=\{m_i,\sigma(m_i)\}.$$
Therefore
$$C_i=\{\ell_i,m_i\}=\{\beta_i,\gamma_i\}=\qtil(B_i)$$
and $\sigma(\ell_i)=m_i, \sigma(m_i)=\ell_i$.
Thus $(Q,\sigma)=(Q_0,\sigma_0)$.
{\enp}
\noindent
Recall Alon's Combinatorial Nullstellensatz (Theorem 1.2 in \cite{Alon99}).
\begin{theorem}[Alon \cite{Alon99}]
\label{t:null}
Let $\FF$ be an be an arbitrary field and let $g=g(x_1,\ldots,x_t) \in \FF[x_1\ldots,x_t]$.
Suppose the total degree $\deg(g)$ of $g$ is $\sum_{i=1}^t d_i$ where each $d_i$ is a nonnegative
integer, and suppose the coefficient of $\prod_{i=1}^t x_i^{d_i}$ in $g$ is nonzero.
Then, if $\Lambda_1,\ldots,\Lambda_t$ are subsets of $\FF$ with $|\Lambda_i|>d_i$, there exist $\lambda_1 \in \Lambda_1,\ldots, \lambda_t \in \Lambda_t$ such that
$g(\lambda_1,\ldots,\lambda_t) \neq 0$.
\end{theorem}
\noindent
{\bf Proof of Theorem \ref{t:rankmatch}:}
By Proposition \ref{p:coeff}, the coefficient of $\prod_{i=1}^t x_i^{\delta_i}$ in $f(x)$ is nonzero.
Applying Theorem \ref{t:null} with $g=f$, $d_i=\delta_i$ and $\Lambda_i=\FF$, and noting that
$\sum_{i=1}^t \delta_i=n=\deg(f)$ and $|\Lambda_i|=|\FF| \geq 3>2 \geq \delta_i$, it follows that there exists a
$\lambda=(\lambda_1,\ldots,\lambda_t) \in \FF^t$ such that $$\det(A+\sum_{i=1}^t \lambda_i B_i)=f(\lambda) \neq 0.$$
{\enp}

\section{Affine Subspaces of Alternating Matrices}
\label{s:alter}
We first recall the definition of the {\it Pfaffian} of an alternating matrix $C=\big(C(i,j)\big)_{i,j=1}^n \in A_{n}(\FF)$ of even order $n=2t$ over an arbitrary field $\FF$.
Let $\ccm_{n}$ denote the set of all perfect matchings in $K_{n}$. For $M=\{e_1,\ldots,e_{t}\} \in \ccm_{n}$, where $e_1 \prec \cdots \prec e_{t}$ and $e_i=\{k_i<\ell_i\}$ for $1 \leq i \leq t$, let
$$\theta(M)=\text{sgn}  \left(
\begin{array}{ccccc}
1 & 2 & \cdots & n-1 & n  \\
k_1 & \ell_1 & \cdots & k_t & \ell_t
\end{array}
\right)
$$
and let
$$\mu(C,M)= \prod_{i=1}^{t} C(k_i,\ell_i).$$
The Pfaffian of $C$ is defined by
$$\pf(C)=\sum_{M \in \ccm_{n}} \theta(M) \mu(C,M).$$
It is well known that $\det(C)=\pf(C)^2$ (see e.g. Exercise 4.24 in \cite{Lovasz}).
\ \\ \\
{\bf Proof of Theorem \ref{t:alter}:} As in the proof of theorem \ref{t:rankmatch}, it suffices to show that if $\ccs$ is an affine subspace of $A_n(\FF)$ and $\nu(G_{\ccs})=\frac{n}{2}=t$, then $\ccs$ contains a nonsingular matrix.
\\
Choose $B_1,\ldots,B_t \in \ccb$ such that
$M_0=\{\qtil(B_1),\ldots,\qtil(B_t)\}$ is a perfect matching of $K_n$.
For $1 \leq i \leq t$ write $q(B_i)=(\alpha_i,\beta_i)$. By reordering we may assume that
\begin{equation}
\label{mono1}
(\alpha_1,\beta_1) \prec  \cdots \prec (\alpha_t,\beta_t).
\end{equation}
Let
$$f(x_1,\ldots,x_t)=\pf(A+\sum_{i=1}^t x_i B_i) \in \FF[x_1,\ldots,x_t].$$
\begin{claim}
\label{c:tfnz}
The coefficient of $x_1\cdots x_t$ in $f(x_1,\ldots,x_t)$ is nonzero.
\end{claim}
\noindent
{\bf Proof:} Let $M=\{e_1,\ldots,e_t\}$ be a matching in $\ccm_{n}$
where $e_i=\{k_i<\ell_i\}$ for $1 \leq i \leq t$ and
\begin{equation}
\label{mono2}
(k_1,\ell_1) \prec \cdots \prec (k_t,\ell_t).
\end{equation}
Then
$$
\mu(A+\sum_{i=1}^t x_i B_i,M)=\prod_{j=1}^t \left(A(k_j,\ell_j)+\sum_{i=1}^t x_i B_i(k_j,\ell_j)\right).
$$
Let $\lambda(M)$ denote the coefficient of $x_1 \cdots x_t$ in $\mu(A+\sum_{i=1}^t x_i B_i,M)$. Then
$$\lambda(M)=\sum_{\pi \in \csy_t} \prod_{j=1}^t B_{\pi(j)}(k_j,\ell_j).$$
Next note that if a permutation $\pi \in \csy_t$ satisfies
$$\prod_{j=1}^t B_{\pi(j)}(k_j,\ell_j) \neq 0$$
then $(k_j,\ell_j) \preceq q(B_{\pi(j)})=(\alpha_{\pi(j)},\beta_{\pi(j)})$ and hence
$$2^{k_j}+2^{\ell_j} \leq 2^{\alpha_{\pi(j)}}+2^{\beta_{\pi(j)}}$$ for all $1 \leq j \leq t$.
Since $$\sum_{j=1}^t \big(2^{k_j}+2^{\ell_j}\big)= \sum_{i=1}^{n} 2^{i}=\sum_{j=1}^t
\big(2^{\alpha_{\pi(j)}}+2^{\beta_{\pi(j)}}\big),$$
it follows that $2^{k_j}+2^{\ell_j} = 2^{\alpha_{\pi(j)}}+2^{\beta_{\pi(j)}}$ and hence
$\{k_j,\ell_j\}=\{\alpha_{\pi(j)},\beta_{\pi(j)}\}$
for all $1 \leq j \leq t$, i.e. $M=M_0$.
The monotonicity assumptions (\ref{mono1}) and (\ref{mono2}) further imply that $\pi$ is the identity permutation.
It follows that the coefficient of $x_1 \cdots x_t$ in $f(x)$ is
$$\theta(M_0) \lambda(M_0)=\theta(M_0) \prod_{j=1}^t  B_j(\alpha_j,\beta_j) \neq 0.$$
{\enp}
\noindent
We now complete the proof of Theorem \ref{t:alter}. By Claim \ref{c:tfnz} the monomial $x_1\cdots x_t$ appears in $f$ with a nonzero coefficient. Applying Theorem \ref{t:null} with $g=f$, $d_i=1$ and $\Lambda_i=\{0,1\}$ for all $i$, and noting that
$\deg(x_1\cdots x_t)=t=\deg(f)$ and $|\Lambda_i|=2 >1=d_i$, it follows that there exists an
$\epsilon=(\epsilon_1,\ldots,\epsilon_t) \in \{0,1\}^t$ such that $f(\epsilon) \neq 0$.
Hence $C=A+\sum_{i=1}^t \epsilon_i B_i \in \ccs$ satisfies
$$\det(C)=\pf(C)^2=f(\epsilon)^2 \neq 0.$$
{\enp}
\noindent
{\bf Proof of Corollary \ref{c:symm}:} Let $\ccs=A+\cb$ where $A\in H_n(\FF_2)$ and $\cb$ is a linear subspace of $H_n(\FF_2)$.
Consider the affine subspace
$$
\ccs'=\Big\{ \left[
\begin{array}{ll}
0 & S \\
S & 0
\end{array}
\right] : S \in \ccs\Big\} \subset A_{2n}(\FF_2).
$$
If $0 \neq B \in \cb$ satisfies $q(B)=(i,j)$ then
$$
q\Big(\left[
\begin{array}{ll}
0 & B \\
B & 0
\end{array}
\right]\Big)= (i,j+n).
$$
It follows that $\nu(G_{\ccs'})=\nu(G_{\ccs})$.
Hence by Theorem \ref{t:alter}
$$
2\rho(\ccs)=\rho(\ccs') \geq \mu(G_{\ccs'})=2\nu(G_{\ccs'})=2\nu(G_{\ccs}).
$$
{\enp}

\section{Maximal Rank and Dimension}
\label{s:eg}
In this section we prove Theorems \ref{t:alt1} and \ref{t:ws1}.
For the alternating case we recall the following extremal graph theoretic result.
\begin{theorem}[Erd\H{o}s-Gallai \cite{EG61}]
\label{t:erga}
Let $G \subset K_n$ be a simple graph that satisfies $\mu(G)=k<n$.
Then
$$
|G| \leq u_a(n,k).
$$
\end{theorem}
\ \\ \\
{\bf Proof of Theorem \ref{t:alt1}:} Let $\ccs \subset A_n(\FF)$ be an affine subspace of $A_n(\FF)$ such that
$\rho(\ccs)=k<n$. Note that $k$ is even. Write $\ccs=A+\ccb$ where $A \in A_n(\FF)$ and $\ccb$ is a linear subspace of $A_n(\FF)$.
By performing Gaussian elimination on a basis of $\ccb$ it is clear that $|G_{\ccs}|=\dim \ccs$.
Theorem \ref{t:alter} implies that
$$
\mu(G_{\ccs}) \leq \rho(\ccs)=k.
$$
It thus follows from Theorem \ref{t:erga} that $\dim \ccs= |G_{\ccs}| \leq u_a(n,k)$.
{\enp}
\noindent
For the weakly symmetric case we need the following version of Theorem \ref{t:erga} for graphs with loops.
\begin{theorem}[\cite{M89}]
\label{t:egm}
Let $G \subset \ktil$ satisfies $\mu(G)=k$, then $|G| \leq u_s(n,k)$.
\end{theorem}
\noindent
{\bf Proof of Theorem \ref{t:ws1}:} Suppose $|\FF|\geq 3$ and
let $\ccs=A+\cb$ where $A \in M_n(\FF)$ and $\cb$ is a linear space contained in $W_n(\FF)$, and let $\rho(\ccs)=k$.
As before $|G_{\ccs}|=\dim \ccs$.
Hence, by Theorem \ref{t:rankmatch}
$$
\mu(G_{\ccs}) \leq \rho(\ccs)=k.
$$
It thus follows from Theorem \ref{t:egm} that $\dim \ccs= |G_{\ccs}| \leq u_s(n,k)$.
{\enp}

\section{Concluding Remarks}
\label{s:conc}

In this note we showed
that if $\ccs$ is a translate of a linear space contained in $W_n(\FF)$ and $|\FF|\geq 3$ then
$\rho(\ccs)\geq \mu(G_{\ccs})$, and that the same holds with no restriction on $\FF$ if $\ccs$ is an affine subspace of $A_n(\FF)$.
These results improve on some earlier work in \cite{M89} and provide a simple approach to upper bounds on the dimension of affine spaces of symmetric and alternating matrices of bounded rank. We conclude with the following two remarks.
\\
1. As noted above, de Seguins Pazzis \cite{DSP16} proved that Theorem \ref{t:ws1} holds also over $\FF_2$ if $\ccs$ is an affine subspace of $H_n(\FF_2)$. It would be interesting to decide whether this case can also be handled using the combinatorial approach of the present paper.
\\
2. When $\ccs$ is a linear subspace of $A_n(\FF)$ and the field $\FF$ satisfies $|\FF|\geq n+1$,
Theorem \ref{t:alter} is the case $p=2$ of Theorem 2.1 in \cite{GM01} that gives a lower bound on the maximal rank of a $p$-vector in a linear subspace $\ccs$ of the $p$-th exterior power $\bigwedge^p \FF^n$ in terms of the weak matching number of a certain $p$-uniform hypergraph associated to $\ccs$. It would be interesting to determine whether this result and its consequences remain true over arbitrary fields.


\begin{thebibliography}{99}

\bibitem{Alon99}
N. Alon,
Combinatorial Nullstellensatz,
{\it Combin. Probab. Comput.}, {\bf 8}(1999) 7–-29.

\bibitem{EG61}
P. Erd\H{o}s and T. Gallai, On the minimal number of vertices representing
the edges of a graph, {\it Publ. Math. Inst. Hungar. Acad. Sci.}, {\bf 6}(1961) 181--203.

\bibitem{GM01}
B. Gelbord and R. Meshulam, Spaces of $p-$vectors of bounded rank,
{\it Israel J. of Math.}, {\bf 121}(2001) 129-139.

\bibitem{Lovasz}
L. Lov\'{a}sz, {\it Combinatorial problems and exercises}, Corrected reprint of the 1993 second edition. AMS Chelsea Publishing, Providence, RI, 2007.

\bibitem{M89}
R. Meshulam, On two extremal matrix problems,
{\it Linear Algebra Appl.}, {\bf 114/115}(1989)  261-271.

\bibitem{DSP16}
C. de Seguins Pazzis,
Affine spaces of symmetric or alternating matrices with bounded rank, {\it Linear Algebra Appl.}, {\bf 504}(2016) 503-–558.

\end{thebibliography}
\end{document}